\documentclass[12pt]{article}

\usepackage[T2A,T1]{fontenc}
\usepackage[utf8]{inputenc}

\usepackage[english,main=russian]{babel}

\usepackage{amsmath, amssymb, amsfonts}
\usepackage{mathtools}

\usepackage[amsmath,amsthm,thmmarks]{ntheorem}

\usepackage[pdftex]{graphicx}
\usepackage{geometry}
\usepackage{enumerate}
\usepackage{colordvi}
\usepackage{ulem}

\usepackage[colorlinks=true,linkcolor=blue,citecolor=blue]{hyperref}

\theoremstyle{plain}
\newtheorem{theorem}{Теорема}[section]      
\newtheorem{lemma}[theorem]{Лемма}          
\newtheorem{corollary}[theorem]{Следствие}

\newtheorem{utv}[theorem]{Утверждение}

\theoremstyle{definition}

\addto\extrasrussian{%
}

\def\F{{\cal F}}
\def\repV{\Upsilon_{\rm v}}
\def\repVSym{\Upsilon_{\rm v}^{\rm sym}}

\let\le\leqslant
\let\leq\leqslant

\usepackage{amsmath}

\begin{document}
\begin{picture}(10.00,00.00)
\put(30.00,40.00){\makebox(0,0)[cc]{УДК 519.157.1+519.175.1+519.176}}
\put(440.00,40.00){\makebox(0,0)[cc]{MSC: 05D15,05C35,05C25}}
\end{picture}
\newline
\begin{minipage}{17cm}
    \begin{center}
        \textbf{\Large{Цена симметрии для хвостатых звёзд}}
    \end{center}
\end{minipage}
\newline
\begin{minipage}{17cm}
    \begin{center}
        Михаил С. Терехов
    \end{center}
\end{minipage}
\newline
\begin{minipage}{17cm}
    \begin{center}
        Механико-математический факультет Московского государственного университета
        Москва 119991, Ленинские Горы, МГУ.
        
        Московский центр фундаментальной и прикладной математики.
    \end{center}
\end{minipage}
\newline
\begin{minipage}{17cm}
    \begin{center}
        kombox.ver.2.0@yandex.ru
    \end{center}
\end{minipage}
\newline

$~~~~~~~~~~$
\begin{minipage}{13.1cm}
\footnotesize{Известно, что если в связном графе $\Gamma$ можно удалить $n$ вершин так, чтобы не осталось подграфов изоморфных графу $K$, 
то можно удалить не более $|V(K)|\cdot n$ вершин, образующих инвариантное относительно всех автоморфизмов графа $\Gamma$ множество, так, чтобы 
не осталось подграфов изоморфных графу $K$. Мы строим бесконечное множество (связных) графов $K$, для которых эта оценка не является точной.
}

\end{minipage}
\newline
\section{Введение}

Рассмотрим следующую задачу.
\newline

\begin{minipage}{15em}

\end{minipage}
$~~~$
\begin{minipage}{15cm}
\normalsize{
	В далёкой галактике решили провести самый большой чемпионат по чтению мыслей. Однако, по регламенту в чемпионате запрещено участвовать 
	галактианцу, который дружит со ста или более участниками. 
	Пусть известно, что как минимум $n$ жителей галактики не смогут участвовать из-за этого правила. 
	Тогда также известно, что если мы хотим сделать множество приглашённых участников инвариантным, относительно автоморфизмов <<графа дружбы>>, то
	получится не приглашать не более $101n$ галактианцев.
}
\end{minipage}
\newline
\newline
\newline
Более общий подход к данной задаче заключён в следующей теореме.
\newline
~~~~~
\newline
\textbf{Теорема KL }[KlLu21]\textbf{.}
\textit{
Пусть группа $G$ действует на множестве $U$ и $\F~-~G$-инвариантное семейство
конечных подмножеств множества $U$, мощности которых ограничены в 
совокупности, а $X$ в $U$ --- конечная система представителей для этого 
семейства (то есть $X \cap F \neq \varnothing$, для всех $F \in \F$). Тогда найдётся 
$G$-инвариантная система представителей Y такая, что $|Y|  \leqslant |X|\max\limits_{F\in \F}|F|.$
}

Слово \textit{семейство} здесь понимается как неупорядоченное
семейство, то есть $\F$ --- это просто некоторое множество подмножеств
множества $U$. Инвариантность семейства $\F$ следует понимать
естественным образом:
${gF\:=\{gf\;|\; f\in F\}\in\F}$ для всех $g\in G$ и $F\in\F$. 

В алгебре есть утверждения, похожие на данную теорему, например, 

-- Если группа $G$ содержит абелеву подгруппу конечного индекса, то $G$ содержит характеристическую (то есть инвариантную относительно всех автоморфизмов) абелеву подгруппу конечного индекса, смотрите, например, [КаМ82].

-- Если группа $G$ содержит разрешимую подгруппу конечного индекса, то $G$ содержит характеристическую разрешимую (той же ступени) подгруппу конечного индекса [KhM07].

-- Если группа $G$ содержит подгруппу конечного индекса с конечным коммутантом, то $G$ содержит характеристическую подгруппу конечного индекса с конечным коммутантом [KlMi15].
\newline

Из теоремы KL легко получается наш случай для графов.
\newline
~~~~~
\newline
\textbf{Следствие KL }[KlLu21]\textbf{.}
\textit{
Пусть $\Gamma$ --- граф и $K$ --- конечный граф. Тогда если в графе $\Gamma$ можно выбрать конечное множество вершин $X$ так,
чтобы каждый подграф графа $\Gamma$, изоморфный графу $K$, имел хоть одну
вершину из $X$, то
в графе $\Gamma$ можно выбрать конечное множество вершин $Y$, инвариантное
относительно всех автоморфизмов графа~$\Gamma$, так, чтобы опять
каждый подграф графа~$\Gamma$, изоморфный графу $K$, имел хоть одну
вершину из $Y$, причём
$|Y| \leqslant |X|\cdot(\hbox{\rm число вершин графа $K$})$.
}
\newline

Здесь слово \textit{граф} мы понимаем как неориентированный граф без петель и кратных рёбер.

В работе [KlLu21] также изучался вопрос о том, является ли оценка из следствия KL неулучшаемой. 
Было доказано, что в общем случае ответ положительный, и поставлен тот же вопрос в классе связных графов. 
В [T22] был получен первый пример графа $K$, для которого оценка не является точной в классе связных графов. 
В данной работе строится бесконечная серия связных графов, для которых указанная оценка не является точной (следствие\autoref{corollary:2.3}).
При этом в процессе доказательства возникла лемма\autoref{lemma:2.1} (о плотных окрестностях), имеющая самостоятельный 
интерес\footnote{А. А. Клячко сообщил автору, что близкий, но менее общий, результат был доказан недавно в [SX25].}.

Мы докажем (доказательство приведено в последнем параграфе) некоторое усиление теоремы KL. Здесь и далее обратный бесконечному кардиналу считается нулём.

\begin{theorem}\label{theorem:1.1}

Пусть группа $G$ действует на множестве $U$ и $\F~-~G$-инвариантное семейство
конечных подмножеств множества $U$, мощности которых ограничены в 
совокупности, а $X$ в $U$ --- конечная система представителей для этого 
семейства (то есть $X \cap F \neq \varnothing$, для всех $F \in \F$). 
Пусть $F\in \F$ имеет непустые пересечения только с орбитами $V_1, ..., V_n$ относительно действия $G$. 
Тогда выполняется неравенство $$\sum_{i=1}^{n}{\frac{|F\cap V_i||X\cap V_i|}{|V_i|}}\geqslant 1.$$
\end{theorem}

Из теоремы\autoref{theorem:1.1} несложно вывести условия на граничные случаи.

\begin{corollary}\label{corollary:1.2}

Пусть группа $G$ действует на множестве $U$ и $\F~-~G$-инвариантное семейство
конечных подмножеств множества $U$, мощности которых ограничены в 
совокупности, а $X$ в $U$ --- минимальная по количеству элементов конечная система представителей для этого 
семейства.
Пусть $X$ имеет непустые пересечения только с орбитами $V_1, ..., V_n$ относительно действия $G$. 
\newline
Тогда

0) $|V_i| < \infty $ для всех $i$.

А если существует минимальное $G$-инвариантное множество $Y$, являющееся системой представителей $\F$ и такое, что 
$|Y|  = |X|\max\limits_{F\in \F}|F|$, то верны следующие утверждения:
\newline

1) Для любой орбиты $V_i$ выполнено $|V_i| = \max\limits_{F\in \F}|F| \cdot |X\cap V_i|$. 

2) $\bigcup\limits_{F\in \F} F \subseteq \bigcup\limits_{i=1}^n V_i$.


3) Для любой орбиты $V_i$ найдётся $F\in \F$, для которого $F \subseteq V_i$.

\end{corollary}
~

Теорема\autoref{theorem:1.1} и следствие\autoref{corollary:1.2} используются для доказательства основного результата (теорема\autoref{theorem:2.4}). 

Выражаю благодарность А.А. Клячко за ценные замечания. Автор благодарит также анонимного рецензента за то, что он внимательно прочитал
более раннюю версию этой статьи и нашёл несколько существенных недочётов. Работа выполнена при поддержке Российского научного фонда, проект № 22-11-00075.
\section{Симметрии на графе}
~
Пусть $\Gamma$ и $K$ --- некоторые графы, 
будем говорить, что множество вершин $M(K) \subseteq V(\Gamma)$ \textit{отмечено} в $\Gamma$, если любой подграф, 
изоморфный графу\footnote{Далее слова ,,изоморфный графу`` будем опускать, если это не вызывает путаницы.} $K$ и 
содержащийся в $\Gamma$, имеет хотя бы одну вершину из $M(K)$. 

Аналогично, будем говорить, что множество вершин $M_{sym}(K) \subseteq V(\Gamma)$ \textit{симметрично отмечено} в $\Gamma$,
если множество $M_{sym}(K)$ является $Aut(\Gamma)$-инвариантным,  
и любой подграф $K$, содержащийся в $\Gamma$, имеет хотя бы одну вершину из $M_{sym}(K)$.

\textit{Вершинной представительностью} $\repV(K,\Gamma)$ графа $K$ в графе $\Gamma$  
называют [KlLu21] минимальное число $n$ такое, что в
графе~$\Gamma$ найдётся отмеченное множество вершин $M(K)$ мощности $n$.

\textit{Симметричной вершинной представительностью} $\repVSym(K,\Gamma)$ графа $K$ в графе $\Gamma$ 
называют [KlLu21] минимальное число $n$ такое, что в графе $\Gamma$ найдётся симметричное отмеченное множество вершин $M_{sym}(K)$ мощности $n$.
\newline

Граф $K$ называют [KlLu21] 
\textit{вершинно-дорогим в классе
графов $\cal K$}, если
$$
\hbox{
$\forall m \in Z$ найдётся граф
$\Gamma_m \in \cal K$ такой, что
}
\infty > \repVSym(K,\Gamma_m)=
\repV(K,\Gamma_m)\cdot\hbox{|V($K$)|} \geqslant m.
$$

Из следствия KL видно, что $\repVSym(K,\Gamma_m) \leq \repV(K,\Gamma_m)\cdot\hbox{|V($K$)|}$.

Граф $K$ называют [KlLu21] \textit{вершинно-дорогим}, если граф $K$ вершинно-дорогой в классе всех графов.

Взяв в качестве $\Gamma_m$ объединение $m$ полных графов $K_{V(K)}$, видно, что любой связный граф $K$ является 
вершинно-дорогим.

\textit{Орбитой}{
$V(v)$ вершины $v$ графа $\Gamma$ назовём $\{gv|g\in Aut(\Gamma)\}$.
}

Граф, индуцированный некоторым подмножеством $V$ из множества вершин графа $\Gamma$, обозначим $\Gamma(V)$.

Будем говорить, что орбита $V$ \textit{пересекается} с $K$, если существует $K \subseteq \Gamma$ такой, что $|V\cap K|>0$.
Далее количество вершин в графе $\Gamma$ будем обозначать $|\Gamma|$.

\begin{utv}\label{utv:2.1}
Пусть $\Gamma$ -- граф, а $K$ -- конечный граф, и для них выполнено условие 
$\repVSym(K,\Gamma)$ $=$ $|K|\cdot\repV(K,\Gamma)$, 
$M(K) \subseteq V(\Gamma)$ -- некоторое (конечное) отмеченное множество в $\Gamma$ с условием $|M(K)| = \repV(K,\Gamma)$. 
Тогда для любой орбиты $V$ графа $\Gamma$ выполнено: если $V$ пересекается с $K$, то 
$\frac{|V \cap M(K)|}{|V|} = \frac{1}{|K|}$, 
иначе $\frac{|V \cap M(K)|}{|V|} = 0$.
\end{utv}
\proof
Если $V$ не пересекается с $K$, то из $M(K)$ можно удалить все вершины $V \cap M(K)$ и множество останется отмеченным, 
откуда из минимальности количества вершин в $M(K)$ следует, что $|V \cap M(K)|=0$. 
Если $V$ пересекается с $K$, то из пункта 2) следствия\autoref{corollary:1.2}, 
применённого к $U=\{$вершины графа $\Gamma\}$, $G = Aut(\Gamma)$, $\F = \{$все подграфы изоморфные $ K \}$ (очевидно, являющееся инвариантным семейством) и $X=M(K)$, получаем,
что $K$ содержится в объединении только тех орбит, в которых есть отмеченные вершины, значит $|V \cap M(K)|>0$, 
откуда из пункта 1) следствия\autoref{corollary:1.2} получаем, что $|V| = |K|\cdot|V \cap M(K)|$.
\endproof

\begin{theorem}\label{theorem:2.2}
\textit{
	Пусть $\Gamma$ -- связный граф, $K$ -- конечный связный граф, имеющий висячую вершину, причём 
	$\infty > |K|\cdot\repV(\Gamma,K) = \repVSym(\Gamma,K)>0$. Тогда для любой орбиты $V$ графа $\Gamma$
	выполнено $\Gamma(V) \supseteq K$.
}
\proof
В самом деле, пусть все орбиты графа $\Gamma$ это $V_{i\in I}, U_{j\in J}$, для некоторых множеств индексов $I,J$, где орбиты $V_i$ пересекаются с $K$, а
орбиты  $U_j$ не пересекаются с $K$. Тогда из пункта 3) следствия\autoref{corollary:1.2} получаем, что все орбиты, которые
пересекаются с $K$, должны содержать $K$. В силу связности графа $\Gamma$ найдутся орбиты $V_i$ и $U_j$, 
между которыми есть ребро в графе $\Gamma$. 
Пусть $K_0$ получается из графа $K$ удалением висячей вершины $x$, соединённой с вершиной $y$.
В силу того, что $V_i$ -- орбита, соединённая каким-то ребром с $U_j$, это ребро можно выбрать так, что один его конец будет в вершине $y$.
Откуда получаем, что $K$ пересекается с $U_j$.
\endproof 
\end{theorem}

\textit{Окрестностью}{
вершины $v$ в графе $\Gamma$ назовём индуцированный подграф графа $\Gamma$, состоящий из всех вершин, 
смежных $v$ и всех рёбер, соединяющих две такие вершины.
}
\newline
\begin{lemma}{\textnormal{(о плотных окрестностях)\textbf{.}}}\label{lemma:2.1}
Не существует непустого графа $\Gamma$, такого, что степень каждой вершины графа равна $n$ и окрестность каждой вершины $v$ графа $\Gamma$ 
содержит $\frac{n(n-1)}{2} - p_v$ рёбер, где $1 \leqslant p_v < \frac{n}{2}$, причём $p_v$ может быть разным для разных вершин.
\proof
Пусть $\Gamma$ -- контрпример. Пусть $v$ -- некоторая вершина, $O(v)$ -- её окрестность. 
Обозначим степень вершины $w$ в подграфе $O(v)$ через $\deg_{O(v)}(w)$. 
Назовём антистепенью вершины $w\in O(v)$ число $\overline{\deg}_{O(v)}(w):=(n-1)-\deg_{O(v)}(w)$. Тогда сумма антистепеней вершин из $O(v)$ 
равна $2p_v$, где $2\leqslant2p_v \leqslant n-1$. Рассмотрим вершину $w$ максимальной антистепени $\overline{\deg}_{O(v)}(w)$ в $O(v)$.
Пусть $\overline{\deg}_{O(v)}(w) = k>0$. Очевидно, что $\overline{\deg}_{O(w)}(v) = \overline{\deg}_{O(v)}(w)$, 
тогда в $O(w)$ есть хотя бы $k$ вершин $w_1, w_2, ..., w_k$: $\overline{\deg}_{O(w)}(w_i)=0$, 
иначе суммарная антистепень вершин из $O(w)$ была бы не меньше, чем $1\cdot k + (k-1)\cdot 0 + (n-1-(k-1))\cdot 1 = n > 2p_w$. 
Вершины $w_1, w_2, ..., w_k$ соединены со всеми вершинами из $O(w)$ (так как $\overline{\deg}_{O(w)}(w_i)=0$), а значит, могут лежать только внутри $O(v) \cap O(w)$ так как $v \in O(w)$, а $v$ соединена только с $O(v)$. 
Тогда вершины $w_1, w_2, ..., w_k$ не соединены с $O(v) \setminus O(w)$ так как все рёбра из них уже ведут в $w$ и $O(w)$. 
Но тогда каждая вершина $u \in O(v) \setminus (O(w)\cup w)$ не соединена с $w, w_1, w_2, ..., w_k$, 
а значит $\overline{\deg}_{O(v)}(u) \geqslant k+1 > \overline{\deg}_{O(v)}(w)$, что противоречит максимальности 
антистепени вершины $w$ в $O(v)$.
\endproof
\end{lemma}
Связный граф c $d+1$ вершиной в котором все рёбра исходят из одной вершины назовём \textit{звездой} и обозначим $S_{d}$.

Граф, получаемый добавлением вершины и ребра, соединяющего её с висячей вершиной графа $S_{d}$ назовём \textit{звездой с хвостом} и обозначим $sr_d$.

\begin{center}
  \setlength{\unitlength}{1cm}

  \begin{minipage}{0.85\textwidth} 

    \begin{picture}(8,3)
      \put(3,1){\circle*{0.1}}  
      \put(5,1){\circle*{0.1}}      
      \put(4,1){\circle*{0.1}}    
      \put(6,1){\circle*{0.1}}  
      \put(5.5,2){\circle*{0.1}} 
      \put(4.5,2){\circle*{0.1}} 
      \put(4.5,0){\circle*{0.1}} 
      \put(5.5,0){\circle*{0.1}}

      \put(3,1){\line(1,0){3}}
      \put(5,1){\line(1,2){0.5}}
      \put(5,1){\line(1,-2){0.5}}
      \put(5,1){\line(-1,2){0.5}}
      \put(5,1){\line(-1,-2){0.5}}

      \put(5,-1){\makebox(0,0)[cc]{$sr_6$}}
    \end{picture}

    \vspace{1cm}

    \setlength{\unitlength}{1mm}
    \linethickness{0.4pt}

    \begin{picture}(20,10)
      \put(85,20){\line(1,1){10}}
      \put(95,30){\line(1,0){20}}
      \put(95,30){\circle*{1}}
      \put(105,30){\circle*{1}}
      \put(115,30){\circle*{1}}
      \put(85,20){\circle*{1}}
      \put(85,40){\line(1,-1){10}}
      \put(85,40){\circle*{1}}
      \put(100,10){\makebox(0,0)[cc]{$sr_3\simeq D_5$}}
      \put(80,5){\makebox(0,0)[cc]{Рис. 1}}
    \end{picture}

  \end{minipage}
\end{center}

Граф на $n$ вершинах со степенью каждой вершины $n-2$ обозначим $K_{n}'$.
\newline
\textbf{Теорема о симметризации систем весовых представителей} [KlT25]\textbf{.}
\textit{Пусть группа $G$ действует на множестве $U$, $W\subset \mathbb{Q}$ — конечное множество (весов), и $\mathcal{F}$ — $G$-инвариантное семейство неотрицательных (весовых) функций $F:U\to W$ с конечным носителем (где $G$-инвариантность означает, что для каждой функции $F\in \F$ и каждого 
$g \in G$ функция $u\to F(g \circ u)$ также лежит в $\F$). Пусть $X\subseteq U$ конечная система весовых представителей для этого семейства, 
то есть $\sum\limits_{x\in X}F(x)\geqslant 1$ для любого $F\in \F$. Тогда найдётся $G$-инвариантная система весовых представителей $Y\subseteq U$
такая, что
$$|Y|\le |X|\cdot \max_{F\in\mathcal{F}} \sum_{u\in U} F(u).$$
При этом в качестве $Y$ можно взять объединение всех $G$-орбит $G\circ u$ таких, что выполнено $|G\circ u \cap X|\cdot \max_{F\in\mathcal{F}} \sum_{u\in U} F(u)  \geqslant |G\circ u|$.
} 
\newline
\begin{utv}\label{utv:2.3}
\textit{
Единственный вершинно-транзитивный связный граф $\Gamma$ со свойством 
\newline
$\infty>|sr_d|\cdot\repV(\Gamma,sr_d) = \repVSym(\Gamma,sr_d)>0$, где $d\geqslant 3$ нечётно --- это граф $K_{d+2}$.
\newline
Существует только 2 вершинно-транзитивных связных графа $\Gamma$ со свойством
\newline
$\infty>|sr_d|\cdot\repV(\Gamma,sr_d) = \repVSym(\Gamma,sr_d)>0$, где $d\geqslant 3$ чётно --- это графы $K_{d+2}$ и $K_{d+2}'$
}
\proof
Пусть $\Gamma$ -- потенциальный контрпример. 
Сначала заметим, что в вершинно-транзитивном графе есть только 2 инвариантных множества вершин -- пустое и все вершины,
значит, $\repVSym(\Gamma,sr_d)\in \{0, |\Gamma|\}$, откуда получаем, что $|\Gamma|<\infty$.
Пусть в графе $\Gamma$ степень каждой вершины равна $k\geqslant d$.
\newline
1) Пусть $k$ нечётно. 
\newline
Тогда по лемме\autoref{lemma:2.1} либо $\Gamma$ является $K_{k+1}$, либо для любой вершины $v$ в $O(v)$ отсутствует хотя бы $\frac{k}{2}$ рёбер.
Тогда в силу нечётности $k$ найдётся вершина $w$ из $O(v)$ которая не соединена хотя бы с двумя вершинами из $O(v)$.
Для вершины $w$ будет выполнено $|O(v) \setminus (O(w) \cup w) |= k-1-|O(v)\cap O(w)| =|O(w) \setminus (O(v) \cup v) | > 1$.
\newline
Зафиксируем две такие вершины $v, u$ и на графе $\Gamma$ введём одну весовую функцию $F$: 
$$
F(x) = \left\{
\begin{aligned}
&1, && \text{если } x \in \{v, w\}, \\
&1, && \text{если } x \in O(v) \cap O(w), \text{но не более } (d-3) \text{ вершин}\\
&\frac{1}{2}, && \text{если }  x \in O(v) \setminus (O(w)\cup w), \text{но не более } (d-1-\min(d-3,|O(v) \cap O(w)|) \text{ вершин}\\
&\frac{1}{2}, && \text{если }  x \in O(w) \setminus (O(v)\cup v), \text{но не более } (d-1-\min(d-3,|O(w) \cap O(v)|) \text{ вершин}\\
& 0,&& \text{для всех остальных вершин}.\\
\end{aligned}
\right.
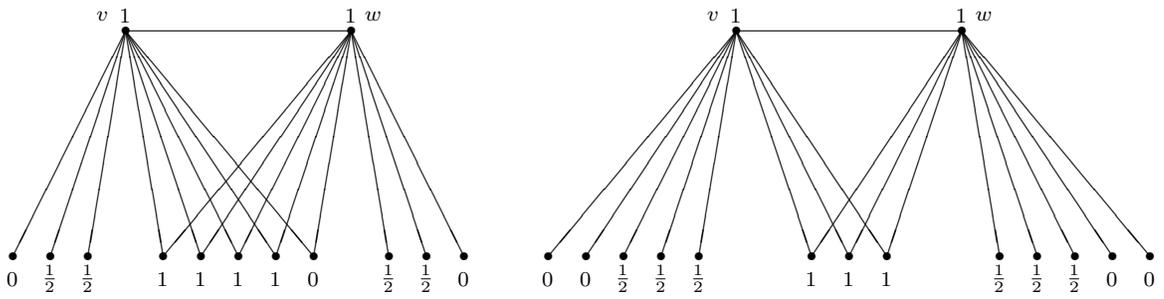
$$
\newline

\unitlength 1.00mm
\linethickness{0.4pt}
\begin{picture}(36.00,60.00)
\put(30.00,48.00){\circle*{1.00}}
\put(60.00,48.00){\circle*{1.00}}
\put(27.00,50.00){\makebox(0,0)[cc]{$_v$}}
\put(63.00,50.00){\makebox(0,0)[cc]{$_w$}}

\put(30.00,50.00){\makebox(0,0)[cc]{$_1$}}
\put(60.00,50.00){\makebox(0,0)[cc]{$_1$}}

\put(30.00,48.00){\line(1,0){30.00}}

\put(30.00,48.00){\line(-1,-6){5.00}}
\put(30.00,48.00){\line(-1,-3){10.00}}
\put(30.00,48.00){\line(-1,-2){15.00}}

\put(30.00,48.00){\line(1,-6){5.00}}
\put(30.00,48.00){\line(1,-3){10.00}}
\put(30.00,48.00){\line(1,-2){15.00}}
\put(30.00,48.00){\line(2,-3){20.00}}
\put(30.00,48.00){\line(5,-6){25.00}}

\put(60.00,48.00){\line(1,-6){5.00}}
\put(60.00,48.00){\line(1,-3){10.00}}
\put(60.00,48.00){\line(1,-2){15.00}}

\put(60.00,48.00){\line(-1,-6){5.00}}
\put(60.00,48.00){\line(-1,-3){10.00}}
\put(60.00,48.00){\line(-1,-2){15.00}}
\put(60.00,48.00){\line(-2,-3){20.00}}
\put(60.00,48.00){\line(-5,-6){25.00}}

\put(25.00,18.00){\circle*{1.00}}
\put(20.00,18.00){\circle*{1.00}}
\put(15.00,18.00){\circle*{1.00}}

\put(35.00,18.00){\circle*{1.00}}
\put(40.00,18.00){\circle*{1.00}}
\put(45.00,18.00){\circle*{1.00}}
\put(50.00,18.00){\circle*{1.00}}
\put(55.00,18.00){\circle*{1.00}}

\put(65.00,18.00){\circle*{1.00}}
\put(70.00,18.00){\circle*{1.00}}
\put(75.00,18.00){\circle*{1.00}}

\put(15.00,15.00){\makebox(0,0)[cc]{$_0$}}
\put(20.00,15.00){\makebox(0,0)[cc]{$_\frac{1}{2}$}}
\put(25.00,15.00){\makebox(0,0)[cc]{$_\frac{1}{2}$}}

\put(65.00,15.00){\makebox(0,0)[cc]{$_\frac{1}{2}$}}
\put(70.00,15.00){\makebox(0,0)[cc]{$_\frac{1}{2}$}}
\put(75.00,15.00){\makebox(0,0)[cc]{$_0$}}

\put(35.00,15.00){\makebox(0,0)[cc]{$_1$}}
\put(40.00,15.00){\makebox(0,0)[cc]{$_1$}}
\put(45.00,15.00){\makebox(0,0)[cc]{$_1$}}
\put(50.00,15.00){\makebox(0,0)[cc]{$_1$}}
\put(55.00,15.00){\makebox(0,0)[cc]{$_0$}}

$~~~~~~~~~~~~~~~~~~~~~~~~~~~~~~~~~~~~~~~~~~~~~~~~~~~~~~~~~~~$
\put(30.00,48.00){\circle*{1.00}}
\put(60.00,48.00){\circle*{1.00}}
\put(27.00,50.00){\makebox(0,0)[cc]{$_v$}}
\put(63.00,50.00){\makebox(0,0)[cc]{$_w$}}

\put(30.00,50.00){\makebox(0,0)[cc]{$_1$}}
\put(60.00,50.00){\makebox(0,0)[cc]{$_1$}}

\put(30.00,48.00){\line(1,0){30.00}}

\put(30.00,48.00){\line(-1,-6){5.00}}
\put(30.00,48.00){\line(-1,-3){10.00}}
\put(30.00,48.00){\line(-1,-2){15.00}}
\put(30.00,48.00){\line(-2,-3){20.00}}
\put(30.00,48.00){\line(-5,-6){25.00}}

\put(30.00,48.00){\line(1,-3){10.00}}
\put(30.00,48.00){\line(1,-2){15.00}}
\put(30.00,48.00){\line(2,-3){20.00}}


\put(60.00,48.00){\line(1,-6){5.00}}
\put(60.00,48.00){\line(1,-3){10.00}}
\put(60.00,48.00){\line(1,-2){15.00}}
\put(60.00,48.00){\line(2,-3){20.00}}
\put(60.00,48.00){\line(5,-6){25.00}}

\put(60.00,48.00){\line(-1,-3){10.00}}
\put(60.00,48.00){\line(-1,-2){15.00}}
\put(60.00,48.00){\line(-2,-3){20.00}}

\put(25.00,18.00){\circle*{1.00}}
\put(20.00,18.00){\circle*{1.00}}
\put(15.00,18.00){\circle*{1.00}}
\put(10.00,18.00){\circle*{1.00}}
\put(5.00,18.00){\circle*{1.00}}

\put(40.00,18.00){\circle*{1.00}}
\put(45.00,18.00){\circle*{1.00}}
\put(50.00,18.00){\circle*{1.00}}

\put(65.00,18.00){\circle*{1.00}}
\put(70.00,18.00){\circle*{1.00}}
\put(75.00,18.00){\circle*{1.00}}
\put(80.00,18.00){\circle*{1.00}}
\put(85.00,18.00){\circle*{1.00}}

\put(5.00,15.00){\makebox(0,0)[cc]{$_0$}}
\put(10.00,15.00){\makebox(0,0)[cc]{$_0$}}
\put(15.00,15.00){\makebox(0,0)[cc]{$_\frac{1}{2}$}}
\put(20.00,15.00){\makebox(0,0)[cc]{$_\frac{1}{2}$}}
\put(25.00,15.00){\makebox(0,0)[cc]{$_\frac{1}{2}$}}

\put(65.00,15.00){\makebox(0,0)[cc]{$_\frac{1}{2}$}}
\put(70.00,15.00){\makebox(0,0)[cc]{$_\frac{1}{2}$}}
\put(75.00,15.00){\makebox(0,0)[cc]{$_\frac{1}{2}$}}
\put(80.00,15.00){\makebox(0,0)[cc]{$_0$}}
\put(85.00,15.00){\makebox(0,0)[cc]{$_0$}}

\put(40.00,15.00){\makebox(0,0)[cc]{$_1$}}
\put(45.00,15.00){\makebox(0,0)[cc]{$_1$}}
\put(50.00,15.00){\makebox(0,0)[cc]{$_1$}}

\put(0.00,5.00){\makebox(0,0)[cc]{$k=9,~d=7$}}
\put(0.00,0.00){\makebox(0,0)[cc]{Рис. 2}}
\end{picture}

Тогда суммарный вес $F$ будет $d+1$, причём любое множество вершин $X$, являющееся системой представителей для
$sr_d$, будет весовой системой представителей для семейства функций $\F:= \{\text{все автоморфные образы } F \}$. 
Действительно, не умаляя общности пусть в $X$ только одна вершина с весом $\frac{1}{2}$, соединённая с $v$. 
Тогда $u$ и все вершины ненулевого веса соединённые с $u$ будут образовывать звезду $S_d$, а в силу того, что $|O(v) \setminus (O(w) \cup w) |>1$, 
найдётся ещё одна вершина веса $\frac{1}{2}$, не принадлежащая $X$ и соединённая с $v$, дополняющая $S_d$ до $sr_d$. 
Тогда по теореме о симметризации систем весовых представителей, существует симметричное множество весовых представителей, не более чем в $d+1$ раз большее, чем 
количество представителей $sr_d$. Однако, единственное симметричное множество в вершинно-транзитивном графе это все вершины (и пустое множество). Значит, $|X|\geqslant \frac{|\Gamma|}{d+1} > \frac{|\Gamma|}{d+2}=\frac{|\Gamma|}{|sr_d|} = \repV(\Gamma,sr_d)$, 
противоречие с существованием отмеченного множества $X'$ со свойством $|X'| = \repV(\Gamma,sr_d)$.
\newline
2) Пусть $k$ чётно. 
\newline
Случай, при котором у каждой вершины $v$ есть 
соединённая с ней вершина $w$ такая, что $|O(v) \setminus (O(w) \cup w) | =|O(w) \setminus (O(v) \cup v) | > 1$, рассматривается аналогично.
Случай, при котором $\Gamma \simeq K_{k+1}$, тоже понятен. Осталось рассмотреть случай, в котором у каждой вершины $v$ окрестность изоморфна $K_k'$. 
В этом случае отсутствующие рёбра разбивают вершины из $O(v)$ на пары. Несложно понять, что это возможно только, когда в графе $k+2$ вершины. 
Причём сам граф -- это граф $K_{k+2}'$. Действительно, рассмотрим 
некоторую вершину $w$, отличную от $v$ и $O(v)$, из которой идёт ребро в вершину $u\in O(v)$. Окрестность $w$ 
должна быть изоморфна $K_k'$, а значит, $u$ должна быть соединена со всеми вершинами из $O(w)$ кроме двух (одна из которых это $u$). 
Однако, все рёбра из $u$ уже известны, а значит, $|O(w)\cap O(v)| \geqslant k-1$. 
Пусть $|O(w)\cap O(v)| = k-1$, но тогда вершина $u' \in O(w)\setminus O(v)$ должна быть соединена с $k-2\geqslant 2$ вершинами из $O(w)$.
Воспользовавшись тем, что у $k-1$ вершины из $O(w)\cap O(v)$ есть рёбра, идущие только в $(O(w)\cap O(v))\cup \{v\} \cup \{u\}$, получаем, 
что $u' \in O(w)\cap O(v)$. Следовательно, $|O(w)\cap O(v)| = k$, а значит, граф -- это граф $K'_{k+2}$.

Но тогда в графе $K'_{k+2}$ требуется отметить хотя бы $k-d+1$ вершину. 
Из условия получаем, $(d+2)\cdot (k-d+1) \leqslant|sr_d|\cdot\repV(\Gamma,sr_d) = \repVSym(\Gamma,sr_d) = k+2$, 
что возможно только в случае $k=d$, а значит, $\Gamma\simeq K'_{k+2} \simeq K_{d+2}'$.
\endproof
\end{utv}
\begin{lemma}\textnormal{[T22]\textbf{.}}\label{lemma:2.5}
Пусть между некоторыми (конечными) орбитами $A$ и $B$ есть хотя бы одно ребро, и 
$ S_1 \subseteq A$ --- некоторое подмножество вершин. Тогда для множества вершин $ S_2 \subseteq B$, связанных хотя бы одним ребром с множеством $S_1$, выполнено неравенство $|S_2| \geqslant |S_1|\cdot\frac{|B|}{|A|}$.
\newline
\end{lemma}

\begin{theorem}\label{theorem:2.4}
Для $n\geqslant 3$ не существует связного графа $\Gamma$ такого, что выполнено:
$$\infty > \repVSym(sr_n,\Gamma)=(n+2)\repV(sr_n,\Gamma)>n+2.$$
\proof
Из утверждения\autoref{utv:2.3} следует, что в графе $\Gamma$ имеется хотя бы две орбиты $A, B$, причём в силу связности их можно выбрать так, что между ними есть хотя бы одно ребро. Из теоремы\autoref{theorem:2.2} следует, что все орбиты содержат $sr_n$, тогда из условия $\infty > \repVSym(sr_n,\Gamma)$, следует, что $|A|, |B|<\infty$. 
Пусть $X$ -- минимальная система представителей для $sr_n$ в графе $\Gamma$, 
тогда по утверждению\autoref{utv:2.1} $|A\cap X| = \frac{|A|}{n+2}$.
Далее получаем цепочку неравенств
$$|A|\stackrel{(1)}{=}\repVSym(sr_n,\Gamma(A))\stackrel{(2)}{\leqslant}|sr_n|\cdot \repV(sr_n,\Gamma(A))\stackrel{(3)}{\leqslant}|sr_n|\cdot |A\cap X| = |A|.$$
Равенство (1) следует из вершинной транзитивности $\Gamma(A)$ и того, что $\Gamma(A)\supseteq s_n.$ 
Неравенство (2) следует из следствия KL.
Неравенство (3) следует из того, что $A\cap X$, являющееся системой представителей для $s_n$ в $\Gamma$, 
является системой представителей для $s_n$ в $\Gamma(A)$.
Тогда неравенство (2) обращается в равенство, и для $\Gamma(A)$ применимо утверждение\autoref{utv:2.3}, откуда следует, 
что каждая орбита является объединением графов $K$ изоморфных $K_{n+2}$ или $K_{n+2}'$. 
Это означает, что для каждого графа $K$ из орбиты $A$ выполнено $|X\cap K|=1$. 
Но тогда существует $|A \setminus X| \geqslant |A|\cdot\frac{n+1}{n+2}$ вершин, 
к которым можно добавить вершину и висячее ребро, и получится граф, содержащий $sr_n$.
Откуда следует, что $X$ содержит все вершины орбиты $B$, связанные хотя бы одним ребром с $A\setminus X$, но таких вершин 
по лемме\autoref{lemma:2.5} будет хотя бы 
$|B|\cdot\frac{n+1}{n+2}$, что невозможно, так как $|B \cap X| = \frac{|B|}{n+2}$.
\endproof
\end{theorem}

\begin{corollary}\label{corollary:2.3}
\textit{
Для $n\geqslant 3$ граф $sr_n$ не является дорогим в классе связных графов.
}
\end{corollary}

\section{Доказательства теоремы 1.1 и следствия 1.2}
Доказательство теоремы\autoref{theorem:1.1} почти дословно повторяет доказательство основной теоремы работы [KlLu21].
\proof
Возьмём какое-то множество $F\in\F$.
Каждое множество $gF$ (где $g\in G$) принадлежит $\F$ в силу
инвариантности семейства $\F$ и, следовательно, пересекается с $X$.
Значит,
$$
G=\bigcup_{f\in F}\{g\in G\;|\;gf\in X\}.
$$
Каждое из множеств $\{g\in G\;|\;gf\in X\}$ является
либо пустым, либо объединением конечного числа левых смежных
классов группы $G$ по стабилизатору $St(f)$ точки $f$:
$$
\{g\in G\;|\;gf\in X\}=
\bigcup_{x\in X}\{g\in G\;|\;gf=x\}=$$

$$=\bigcup_{x\in X\cap Gf}g_x\cdot St(f),~
\{\text{где } g_x\in G \text{ фиксированы так, что } g_xf=x\}.
\eqno{(1)}
$$
Таким образом, мы получили разложение группы~$G$ в конечное объединение
левых смежных классов по некоторым подгруппам.
Воспользуемся теперь теоремой Б.~Неймана [Neu54] (утверждение~4.5):
если группа $G$ покрывается конечным числом смежных классов по некоторым
\(необязательно разным\) подгруппам: $G=g_1G_1\cup\dots\cup g_sG_s$,
то
\
$\displaystyle\sum {1\over|G:G_i|} \geqslant 1$
\(где обратный к бесконечному кардиналу
считается нулём\).
Следовательно, (учитывая то, что индекс стабилизатора равен длине орбиты)
мы получаем 
$$1 \leqslant \sum_{f\in F}{1\over|G: St(f)|}\cdot|Gf\cap X|=
\sum_{f\in F}{|Gf\cap X|\over|Gf|} = 
\sum_{i=1}^{n}\sum_{f\in F\cap V_i}{|Gf\cap X|\over|Gf|} = 
\sum_{i=1}^{n}{\frac{|F\cap V_i||X\cap V_i|}{|V_i|}}.\quad \fbox{}$$

Теперь докажем следствие\autoref{corollary:1.2}.
\proof
$~$

0) Возьмём какое-то множество $F\in\F$.
Пусть $F = F_1\sqcup F_2$, 
где $F_1$ -- элементы, принадлежащие конечным орбитам, а $F_2$ -- бесконечным. Аналогично (1) из теоремы 1.1 получаем 
разложение группы~$G$ в конечное объединение левых смежных классов по некоторым подгруппам.
Воспользуемся теперь теоремой Б.~Неймана [Neu54] (утверждение~4.4):
Если группа покрыта конечным числом смежных классов по подгруппам,
то все смежные классы по подгруппам бесконечного индекса
можно из этого покрытия выбросить, и оно всё равно останется покрытием.  
$$G = \bigcup_{f\in F}\bigcup_{x\in X\cap Gf}g_x\cdot St(f) = \bigcup_{f\in F_1}\bigcup_{x\in X\cap Gf}g_x\cdot St(f) = $$
$$= \bigcup_{f\in F_1}\bigcup_{x\in X}\{g\in G\;|\;gf=x\} = \bigcup_{f\in F_1}\{g\in G\;|\;gf\in X\}. $$
Из принадлежности единицы группы $G$ последнему множеству получаем, что $f = ef \in X$ для некоторого $f\in F_1$. 
Откуда следует, что $|X\cap F_1|>0$, а это значит, что отбросив элементы $X$ из бесконечных 
орбит мы, по-прежнему получаем систему представителей. 

1) Пусть $M := \max\limits_{F\in \F}|F|$.
Из теоремы о симметризации систем весовых представителей найдется симметричная система представителей $Y'$, такая, 
что для каждой орбиты $V$ выполнено $|Y'\cap V| \leqslant M\cdot|X\cap V|$. Тогда $M\cdot |X| = |Y|\leqslant |Y'| \leqslant M\cdot |X|$, 
откуда следует, что все неравенства обращаются в равенство.

2) Рассмотрим некоторый $F \in \F$ и подставим равенство из пункта 1) в теорему\autoref{theorem:1.1} 
$$1 \leqslant \sum_{i=1}^{n}{\frac{|F\cap V_i||X\cap V_i|}{|V_i|}} = \sum_{i=1}^{n}{\frac{|F\cap V_i|}{M}} \leqslant 1.$$
Откуда следует утверждение пункта 2).

3) Пусть орбита $V_1$ не содержит никакой $F \in \F$. Тогда из утверждения пункта 2) легко видеть, что 
$Y' := \bigcup\limits_{i=2}^n V_i$ является симметричной системой представителей. Применяя пункт 1), получаем  
$|Y| \leqslant |Y'| < |\bigcup\limits_{i=1}^n V_i| = M \cdot |\bigcup\limits_{i=1}^n (X\cap V_i)| = M\cdot |X|$. Противоречие.
\endproof

СПИСОК ЦИТИРОВАННОЙ ЛИТЕРАТУРЫ
\newline
[KlLu21] A. A. Klyachko, N. M. Luneva, Invariant systems of representatives, or The cost of symmetry,
Discrete Mathematics, 344:6 (2021), 112361. См. также arXiv:1908.03315.
\newline
[КаМ82] М. И. Карграполов, Ю. И. Мерзляков, Основы теории групп. М.: Наука, 1982.
\newline
[KhM07] E. I. Khukhro, N. Yu. Makarenko,
Large characteristic subgroups satisfying multilinear commutator
identities,
J. London Math. Soc., 75:3 (2007), 635-646.
\newline
[KlMi15] A. A. Klyachko, M. V. Milentyeva,
Large and symmetric: The Khukhro-Makarenko theorem on laws --- without laws,
J. Algebra, 424 (2015), 222-241. arXiv 1309.0571
\newline
[KlT25]  A. A. Klyachko, M. S. Terekhov, Invariant systems of weighted representatives. J Algebr Comb 61, 32 (2025). https://doi.org/10.1007/s10801-025-01400-y
\newline
[Neu54] B.H. Neumann, Groups covered by 
permutable subsets, J. London Math. Soc., s1-29:2 (1954), 236-248.
\newline
[SX25] N. S. Sheffield, Z. Xi, 
Graphs with the same edge count in each neighborhood, 
\newline
arxiv.org/abs/2507.14473, 2025.
\newline
[T22] М.С. Терехов, Цена симметрии в связных графах, Мат. заметки,
112:6 (2022), 895-902. См. также arXiv:2202.09590
\end{document}